\documentclass[aop,MSNbibl,dvips]{arximspdf}

\doi{10.1214/09-AOP515}
\volume{38}
\issue{4}
\pubyear{2010}
\firstpage{1390}
\lastpage{1400}

\makeatletter

\newtheorem{theorem}{Theorem}
\newproclaim{remark}[theorem]{Remark}
\newtheorem{corollary}[theorem]{Corollary}
\newtheorem{lemma}[theorem]{Lemma}
\newtheorem{proposition}[theorem]{Proposition}

\makeatother

\begin{document}
\begin{frontmatter}

\title{Right inverses of L\'{e}vy processes}
\runtitle{L\'{e}vy inverses}

\pdftitle{Right inverses of Levy processes}

\begin{aug}
\author[A]{\fnms{Ron} \snm{Doney}\corref{}\ead[label=e1]{rad@ma.man.ac.uk}} and
\author[B]{\fnms{Mladen} \snm{Savov}\ead[label=e2]{savov@stats.ox.ac.uk}}
\runauthor{R. Doney and M. Savov}
\affiliation{University of Manchester and Oxford University}
\address[A]{School of Mathematics\\
University of Manchester\\
Oxford Road\\
Manchester M13 9PL\\
United Kingdom\\
\printead{e1}} 
\address[B]{Department of Statistics\\
Oxford University\\
1 South Parks Road\\
Oxford, OX2 3TG\\
United Kingdom\\
\printead{e2}}
\end{aug}

\received{\smonth{4} \syear{2009}}
\revised{\smonth{11} \syear{2009}}

%
\begin{abstract}
We call a right-continuous increasing process $K_{x}$ a \textit{partial right
inverse} (PRI) of a given L\'{e}vy process $X$ if $X_{K_{x}}=x$ for at least
all $x$ in some random interval $[0,\zeta)$ of positive length. In this
paper, we give a necessary and sufficient condition for the existence
of a
PRI in terms of the L\'{e}vy triplet.
\end{abstract}

\setattribute{keyword}{AMS}{AMS 2000 subject classification.}

\pdfkeywords{60G51, Levy process, ladder height subordinator,
sample path behavior, creeping, excursions}

\end{frontmatter}

\section{Introduction and results}

In this paper, a real-valued L\'{e}vy process is studied. The problem of
existence of a partial right inverse (PRI) is considered and an explicit
integral criterion is provided for testing whether any L\'{e}vy process
possesses a PRI.

We continue work by Evans \cite{se} and Winkel \cite{mw}. Evans has
introduced the notion of a full right inverse and has defined this
process $K$
as the minimal increasing process that satisfies $X(K_{x})=x$ for all $%
x\geq0$; Winkel, in \cite{mw}, has extended this definition to $%
X(K_{x})=x $ on some random interval $[0,\zeta)$ of positive length and
has named this process a PRI. In these two papers, it is shown that if $K$
exists, it is a (possibly killed) subordinator.

A L\'{e}vy process $X=(X_{t};t\geq0)$ is a stochastic process which
possesses stationary and independent increments, starts from zero and whose
paths are a.s. right-continuous. Each L\'{e}vy process is fully
characterized by its L\'{e}vy triplet $(\gamma,\sigma,\Pi)$, where $%
\gamma\in\mathbb{R}$, $\sigma\geq0$ and the L\'{e}vy measure $\Pi$ has
the property
\[
\int_{-\infty}^{\infty}(1\wedge x^{2})\Pi(dx)<\infty.
\]
Also, each L\'{e}vy process $X$ can be represented as follows:
%
%
\begin{equation} \label{dec}
X_{t}=\gamma t+\sigma B_{t}+X_{t}^{(1)}+\sum_{0<s\leq t}(X_{s}-X_{s-})%
\mathbf{1}_{(|X_{s}-X_{s-}|>1)},
\end{equation}
where $B$ is a standard Brownian motion, $X^{(1)}$ is a pure jump zero mean
martingale and all of the components in (\ref{dec}) are independent. In the
class of L\'{e}vy processes, we distinguish between L\'{e}vy processes with
bounded variation and L\'{e}vy processes with unbounded variation. The former
are those for which $\sigma=0$ and $\int_{-\infty}^{\infty}(1\wedge
|x|)\Pi(dx)<\infty$. In this case, $X$ can be represented as
%
%
\begin{equation} \label{bv}
X_{t}=bt+X_{t}^{+}+X_{t}^{-},
\end{equation}
where $b$ is the drift coefficient and $X^{+}$ and $X^{-}$ are independent
driftless subordinators (i.e., increasing L\'{e}vy processes). In our setting,
as well as in many other situations, these two classes of processes exhibit
quite different behaviors and need separate attention.

We write $R_{t}=\sup_{s\leq t}X_{s}-X_{t}$. It is shown in \cite{jb}, Chapter
6, that $R$ is a strong Markov process which possesses a local time at
zero, $L(t)$, and a corresponding inverse local time $L^{-1}(t)=\inf
\{s\dvtx L(s)>t\}$ such that $(L^{-1}(t);X(L^{-1}(t)))$ is a bivariate
subordinator: we denote its L\'{e}vy measure by $\mu^{(+)}(dt;dy)$ and we
use, in particular, $\mu^{(+)}(dy)=\mu^{(+)}((0;\infty);dy)$. We
also use
the notation $H^{+}(t):=X(L^{-1}(t))$ and call $H^{+}$ the \textit{upward ladder
height process}. Similarly, we can define $Z_{t}=X_{t}-\inf_{s\leq t}X_{s}$
and, using the same arguments, we have an associated inverse local time
$%
L_{-}^{-1}(t)$ and \textit{downward ladder height process} $%
H^{-}(t):=X(L_{-}^{-1}(t))$. We denote the L\'{e}vy measure of $H^{-}$
by $%
\mu^{(-)}(dy)$. Finally, with each of the subordinators $H^{+}$and
$H^{-}$, we
associate the so-called \textit{renewal measure}, defined as follows:
%
%
\begin{equation}\label{ren}
U_{+}(x)=E\int_{0}^{\infty}\mathbf{1}_{\{H_{t}^{+}\leq x\}}\,dt,\qquad
U_{-}(x)=E\int_{0}^{\infty}\mathbf{1}_{\{H_{t}^{-}\leq x\}}\,dt.
\end{equation}
We refer to Bertoin \cite{jb} or Doney \cite{rad} for more information
on L%
\'{e}vy processes.

Next, we briefly discuss the definition of a PRI, that is,
$K=(K_{x},x\geq0)$.
We follow an approach developed in Evans \cite{se}. Define, for each
$n\geq
1 $ and $k\geq0$, the stopping times
%
%
\begin{equation}\label{k}
T_{0}=0,\qquad T_{n}^{k+1}=\inf\biggl\{t\geq T_{n}^{k}\dvtx X_{t}=\frac{k+1}{2^{n}}\biggr\}
\end{equation}
and processes
\[
K_{x}^{n}=T_{n}^{k},\qquad \frac{k}{2^{n}}\leq x<\frac{k+1}{2^{n}}.
\]
A pathwise argument then shows that
%
%
\begin{equation} \label{ev}
K_{x}=\inf_{y>x}\sup_{n\geq0}K_{y}^{n}.
\end{equation}
It is possible that for each $x>0$, the definition above gives $K_{x}%
\stackrel{\mathrm{a.s.}}{=}\infty$ and, in this case, we say that a
PRI does not
exist. The question of the existence of a PRI has been studied by Evans
in \cite%
{se} and Winkel in \cite{mw}. Evans has shown that for any symmetric L\'
{e}%
vy process with $\sigma>0$, a full right-inverse exists. Winkel
\cite{mw} then showed that the same result holds for any oscillating L\'
{e}vy process
with $\sigma>0$ and also described all L\'{e}vy processes with bounded
variation having a PRI. Moreover, in the unbounded variation case, he
provided a necessary and sufficient condition (NASC) for the existence
of a
PRI, but this NASC is not satisfactory since it requires knowledge
about the
second derivative at zero of the so-called \mbox{$q$-potentials} of the given L\'
{e}%
vy process, which are generally unknown. Therefore, the main aim of this
paper is to supply an NASC for the existence of a PRI in terms of the
L\'{e}vy
triplet, that is, $(\gamma,\sigma,\Pi)$, in the unbounded variation
case.
In fact, our method, which is probabilistic in nature, also deals with the
bounded variation case and gives the following result.
\begin{theorem}
\label{T}
Let $X$ be a L\'{e}vy process with a L\'{e}vy measure $\Pi$ such
that $\Pi(\mathbb{R})>0$. Then:

\begin{longlist}
\item if $X$ has unbounded variation, it has a partial right inverse
(PRI) iff
\newline
$\sigma>0$ or $\sigma=0$, $\Pi(\mathbb{R}^{-})=\infty$ and
$J<\infty$, where, with
$\overline{\Pi}{}^{(-)}(s)=\int_{-\infty}^{-s}\Pi(dx)$,
%
%
\begin{equation} \label{J}
J=\int_{0}^{1}\frac{x^{2}\Pi(dx)}{ ( \int_{0}^{x}\int
_{y}^{1}\overline{%
\Pi}{}^{(-)}(s)\,ds\,dy ) ^{2}} ;
\end{equation}

\item if $X$ has bounded variation, then it has a PRI iff $\Pi
(\mathbb{R}%
^{+})<\infty$ and $X$ has a drift coefficient $b>0$.
\end{longlist}
\end{theorem}
\begin{remark}
If $\Pi(\mathbb{R)}=0$, then $X_{t}=\gamma t+\sigma B_{t}$ is a continuous
process and $T_{x}=\inf\{t\dvtx X_{t}=x\}$ will be a PRI on the set $%
\{T_{x}<\infty\}$. Note that, in this case, $\{T_{x}<\infty\}$ will
be the
empty set iff $\sigma=0$ and $\gamma<0$.
\end{remark}
\begin{remark}
A L\'{e}vy process $X$ is said to ``creep upward'' if $P(X(T_{x}^{+})=x)>0$
for some (and then all) $x>0$, where $T_{x}^{+}=\inf(t>0\dvtx X_{t}>x)$. It is
known that this happens iff the ladder height process $H^{+}$ has drift
$%
\delta_{+}>0$; see, for example, Theorem 19, page 174 of \cite{jb}.
Since it is always the case that $%
\sigma^{2}=2\delta_{+}\delta_{-}$, where $\delta_{-}$ is the drift
of $%
H^{-}$, this certainly happens when $\sigma>0$. If $\sigma=0$ and $%
J<\infty$, then the integral
%
%
\begin{equation}\label{K}
L=\int_{0}^{1}\frac{x^{2}\Pi(dx)}{\int_{0}^{x}\int_{y}^{1}\overline
{\Pi}{}
^{(-)}(s)\,ds\,dy}
\end{equation}
is clearly finite and it is shown in \cite{vv} that this is the NASC
for $\delta
_{+}>0$ in the unbounded variation case when $\sigma=0$. (See also Section
6.4 of \cite{rad} for an alternative proof of this result.) Finally, in the
bounded variation case, $b>0$ is clearly equivalent to $\delta_{+}>0$.
We therefore
conclude that our theorem is consistent with the intuitively obvious claim
that ``upward creeping'' is necessary, but not sufficient, for the existence
of a PRI.
\end{remark}

The next corollary illustrates how our theorem yields specific information
in special cases. Here, and throughout the paper, we use the notation $%
f\approx g$ to denote the existence of constants $0<c<C<\infty$ with $%
cg(x)\leq f(x)\leq Cg(x)$, for all sufficiently small $x$.
\begin{corollary}
\label{Co}Let $X$ be a L\'{e}vy process with $\sigma=0$ and L\'{e}vy
measure $\Pi$ such that $\overline{\Pi}{}^{+}(x)=\int_{x}^{\infty}\Pi
(dy)\approx x^{-\beta}$ and $\overline{\Pi}{}^{-}(x)\approx x^{-\alpha}$,
where $1\leq\alpha<2$ and $0\leq\beta<2$. Then $X$ has a PRI iff
$\beta
<2\alpha-2$.
\end{corollary}
\begin{remark}
This result extends Proposition 2 and Theorem 6 in \cite{mw}.
\end{remark}

\section{Proofs}

Recall that we denote by $H^{+}$ the ascending ladder height process of a
given L\'{e}vy process $X$. We use $\delta_{+}$ to denote the drift of
$%
H^{+}$ and $\mu^{(+)}(dy)$ to denote its L\'{e}vy measure. We also use
$%
U_{+}$ and $U_{-}$, which are defined in (\ref{ren}). We start the
proof by
disposing of some special cases.

Suppose, first, that $\Pi(\mathbb{R})<\infty$. Then $V=\inf
\{t>0\dvtx X_{t}-X_{t-}\neq0\}>0$ a.s. since it is an exponentially distributed
random variable with parameter $\Pi(\mathbb{R})$ and the given process
coincides up to time $V$ with the process we get by removing all of its jumps.
The resulting process will be of the form $\sigma B_{t}+bt$, which possesses
a PRI iff $\sigma>0$ or $\sigma=0$ and $b>0$, in accordance with
Theorem \ref{T}. Next, suppose that $\Pi(\mathbb{R})=\infty$, but $\Pi(\mathbb
{R}%
^{+})<\infty$. Removing all the positive jumps then gives a spectrally
negative L\'{e}vy process $\tilde{X}$. If $\tilde{X}$ has unbounded
variation, or has bounded variation and a positive drift $b$, then it passes
continuously over positive levels. Then\vspace*{-1pt} with $\tilde{T}(x)=\inf\{
t>0\dvtx
\tilde{X}_{t}=x\}$, we obviously have $\tilde{X}_{\tilde{T}(x)}=x$ on
$\{%
\tilde{T}(x)<\infty\}$ and we can choose $K_{x}=\tilde{T}(x)$.
Alternatively, $\tilde{X}$ has bounded variation and a drift $b\leq0$, and,
clearly, no PRI exists for $\tilde{X}$ or $X$ in this case. Noting that in
the unbounded variation case, we have $\int_{0}^{1}\overline{\Pi}{}%
^{(-)}(s)\,ds=\infty$ so that, necessarily, $J<\infty$, we see that these
results also accord with Theorem \ref{T}. Next, suppose that $\Pi
(\mathbb{R}%
)=\infty$, but $\Pi(\mathbb{R}^{-})<\infty$. If $X$ has bounded
variation, then removing all of the negative jumps gives us a spectrally
positive process of the form $\tilde{X}_{t}=X_{t}^{+}+bt$, where
$X^{+}$ is
a driftless subordinator. If $b\geq0$, then $\tilde{X}$ has monotone paths
and the assumption that $\Pi(\mathbb{R}^{+})=\infty$ implies the
existence of points $x_{n}\downarrow0$ with $P(T(x_{n})=\infty)=1$, which
verifies Theorem \ref{T} in this case. Finally, if $b<0$ or if $X$ has
unbounded variation, then the decreasing ladder height process is a pure
drift, possibly killed at an exponential time, and we see that the
hypothesis of Proposition \ref{A} below holds.

The rest of our proof uses the following simple consequence of the
construction of $K$ due to Evans \cite{se}.
\begin{lemma}
\label{O}Let $X$ be an arbitrary L\'{e}vy process, and set $T_{x}=\inf
\{t>0\dvtx X_{t}=x\}$ and $p_{x}=P(T_{x}=\infty)=P(X$ does not visit $x)$. Then:

\begin{longlist}
\item a PRI exists for $X$ if%
%
%
\begin{equation} \label{b}
\lim\sup_{x\downarrow0}\frac{1-E(e^{-\theta T_{x}})}{x}<\infty\qquad\mbox{for
some }\theta>0;
\end{equation}

\item no PRI exists for $X$ if
%
%
\begin{equation} \label{c}
\lim_{x\downarrow0}x^{-1}p_{x}=\infty.
\end{equation}
\end{longlist}
\end{lemma}
\begin{pf}
First, note that the sequence $K^{(n)}:=T_{n}^{2^{n}}$, $n\geq1$,
where $%
T_{n}^{k}$ are defined in (\ref{k}), is monotone increasing. If we denote
its limit by $\tilde{K}$, then it is immediate from (\ref{ev}) that $%
K_{1}\leq\tilde{K}\leq K_{2}$. Since we know that $K$ is a (possibly
killed) subordinator, we see that existence of a PRI for $X$ is equivalent
to $P(\tilde{K}<\infty)>0$. However, this is equivalent to
\[
\lim_{n\rightarrow\infty}E\bigl(e^{-\theta K^{(n)}}\bigr)=E(e^{-\theta\tilde
{K}}\dvtx
\tilde{K}<\infty)>0
\]
for some (and then all) $\theta>0$. Since $K^{(n)}$ is the sum of $2^{n}$
independent random variables distributed as $T_{2^{-n}}$, we have
\[
\log E(e^{-\theta\tilde{K}}\dvtx\tilde{K}<\infty)=\lim_{n\rightarrow
\infty
}2^{n}\log E(e^{-\theta T_{2^{-n}}})
\]
and this is clearly finite for any $\theta$ for which (\ref{b}) holds.
Since $1-E(e^{-\theta T_{x}})\geq p_{x}$, we see that this limit is
$-\infty$
for all $\theta>0$ whenever (\ref{c}) holds, and the result follows.
\end{pf}

The crux of our proof is contained in the following result, in which $%
\overline{\mu}{}^{+}(x)=\mu((x,\infty))$ for $x>0$.
\begin{proposition}
\label{A} Let $X$ be a L\'{e}vy process having $\Pi(\mathbb
{R}^{+})=\infty$
and $U_{-}(dx)>0$ for all small enough $x>0$. Then $X$ has a PRI iff
$\delta
_{+}>0$ and $I<\infty$, where
%
%
\begin{equation} \label{I}
I=\int_{0}^{1}\overline{\mu}{}^{+}(x)U_{-}(dx)=\int_{0}^{1}\mu
^{+}(dx)U_{-}(x).
\end{equation}
\end{proposition}
\begin{pf}
Since the existence of a PRI is a local property, we can truncate the
L\'{e}vy
measure so that it is contained in $[-1;1]$. Indeed, the first jump of $X$
larger than $1$ in absolute value occurs after an exponential time
$\zeta$
and $K_{x}$ is a subordinator, therefore $K_{x}<\zeta$ pathwise for
all $x$ small enough. This shows that the existence of $K$ is
independent of
the large jumps, so we will assume, without loss of generality, that
$\Pi([1,\infty))=\Pi
((-\infty,-1]))=0$. Moreover, the value of $\delta_{+}$ is also a
local property, so this is also unchanged by any alteration of the L\'{e}vy
measure on closed intervals which do not contain 0. Note that our
assumptions imply that $I>0$ and that these alterations do not change the
finiteness/infiniteness of $I$. Let us introduce some notation. For
$x>0$, we
put $T_{x}^{+}=\inf\{t>0\dvtx X_{t}>x\}$ and $T_{x}^{-}=\inf\{
t>0\dvtx X_{t}<-x\}$
for the first passage times above $x$ and below $-x$, respectively, and
$%
O^{+}(x)=X_{T_{x}^{+}}-x$, $O^{-}(x)=x-X_{T_{x}^{-}}$ for the overshoot
above $x$ and the undershoot below $-x$, respectively. Noting that
$O^{+}(x)$ is also the
overshoot of $H^{+}$ above $x$, we can use Proposition 2, page 76 in
\cite{jb}
to deduce that for $x>0$, $y>0$,%
%
%
\begin{eqnarray} \label{bound}
\overline{\mu}{}^{(+)}(x+y)U_{+}(x) &\leq&P\bigl(O^{+}(x)>y\bigr)=\int_{0}^{x}%
\overline{\mu}^{(+)}(x+y-z)U_{+}(dz) \nonumber\\[-8pt]\\[-8pt]
&\leq&\overline{\mu}^{(+)}(y)U_{+}(x).\nonumber
\end{eqnarray}
To prove the result in one direction, we alter the L\'{e}vy measure by
adding a mass at $\{1\}$, if necessary, to make $X$ drift to $+\infty
$. We
then have the estimate%
\begin{eqnarray*}
p_{x} &\geq&P(O_{x}^{+}>0,\mbox{ and }X\mbox{ stays above }x) \\
&=&\int_{0}^{1}P\bigl(O^{+}(x)\in dy\bigr)P(T_{y}^{-}=\infty) \\
&=&c\int_{0}^{1}P\bigl(O^{+}(x)\in dy\bigr)U_{-}(y) \\
&=&c\int_{0}^{1}P\bigl(O^{+}(x)>y\bigr)U_{-}(dy),
\end{eqnarray*}
where the fact that $P(T_{y}^{-}=\infty)=cU_{-}(y)$ comes from Proposition
17, page 172 of \cite{jb}. [It is obvious that, in fact,
$c=1/U_{-}(\infty)$
since $P(T_{y}^{-}=\infty)\rightarrow1$ as $y\rightarrow\infty$.]
From (\ref{bound}), it then follows that%
\begin{eqnarray*}
\lim_{x\downarrow0}\inf x^{-1}p_{x} &\geq&c\lim_{x\downarrow0}\inf
x^{-1}U_{+}(x)\int_{0}^{1}\overline{\mu}^{(+)}(x+y)U_{-}(dy) \\
&\geq&cI\lim_{x\downarrow0}\inf x^{-1}U_{+}(x).
\end{eqnarray*}
Finally, we recall from Proposition 1, page 74 in \cite{jb} that $%
U_{+}(x)\approx x/(\delta_{+}+\int_{0}^{x}\overline{\mu}^{(+)}(y)\,dy)$
so that $x^{-1}U_{+}(x)\approx1/\delta_{+}$ as $x\downarrow0$, and
thus (%
\ref{c}) holds and no PRI exists, whenever $\delta_{+}=0$, or $\delta
_{+}>0 $ and $I=\infty$. To argue in the other direction, we assume
that $%
\delta_{+}>0$ and $I<\infty$. Then, without loss of generality, we
can take $%
\delta_{+}=1$. Next, we denote by $P^{\theta}$ the law of this process
killed at an independent exponential time $\tau$ with parameter
$\theta$
and note that
\[
p_{x}^{\theta}:=P^{\theta}(T_{x}=\infty)=P(T_{x}>\tau
)=1-E(e^{-\theta
T_{x}}).
\]
Our aim is to show that there exists some $ \theta>0$ such that
%
%
\begin{equation} \label{lim}
\lim\sup_{x\downarrow0}x^{-1}p_{x}^{\theta}<\infty
\end{equation}
since then the existence of a PRI for $X$ will follow from Lemma \ref
{O}. We
decompose $p_{x}^{\theta}$ according to the number of upcrossings and
downcrossings of level $x$ that occur. To do so, we denote by $T^{+}(x,n)$
the time of $n$th crossing above $x$, by $T^{-}(x,n)$ the time of $n$th
crossing below $x$ and for $n\geq1$, we put
\begin{eqnarray*}
p_{x}^{\theta}(n) &=&P^{\theta}\bigl\{T_{x}=\infty,T^{+}(x,n)<\infty
,T^{-}(x,n)=\infty\bigr) ,\\
q_{x}^{\theta}(n) &=&P^{\theta}\bigl\{T_{x}=\infty,T^{-}(x,n)<\infty
,T^{+}(x,n+1)=\infty\bigr).
\end{eqnarray*}
Since $X$ creeps upward, it is then easy to see that%
%
%
\begin{equation} \label{sum}
p_{x}^{\theta}=P^{\theta}(T_{x}^{+}=\infty)+\sum_{1}^{\infty
}p_{x}^{\theta}(n)+\sum_{1}^{\infty}q_{x}^{\theta}(n).
\end{equation}
We start by noting that
\[
P^{\theta}(T_{x}^{+}=\infty)=c^{+}(\theta)U_{+}^{\theta}(x)\qquad\mbox{where
}c^{+}(\theta)=\frac{1}{U_{+}^{\theta}(\infty)},
\]
and $U_{+}^{\theta}(x)$ is the renewal function of the ladder height
process $H^{+}$ under $P^{\theta}$. Of course, under $P^{\theta}$, $H^{+}$
is killed at some rate $k^{+}(\theta)>0$ and has L\'{e}vy measure $\mu
^{+}(\theta,dx)\leq\mu^{+}(dx)$. However, as we have mentioned, its
drift is
unchanged and equals $1$. Using a version of Erickson's bound for killed
subordinators, which can be found in \cite{vv}, we therefore have
%
%
\begin{equation} \label{er}
U_{+}^{\theta}(x)\leq\frac{c_{0}x}{1+\int_{0}^{x}\overline{\mu}%
^{+}(y,\theta)\,dy+xk^{+}(\theta)}\leq c_{0}x,
\end{equation}
where $c_{0}$ is an absolute constant. Also,
\[
U_{+}^{\theta}(\infty)=\lim_{y\rightarrow\infty}\int_{0}^{\infty
}e^{-tk^{+}(\theta)}P(H_{t}^{+}\leq y)\,dt=\frac{1}{k^{+}(\theta)}
\]
and this gives the bound%
%
%
\begin{equation} \label{zero}
P^{\theta}(T_{x}^{+}=\infty)\leq c_{0}k^{+}(\theta)x.
\end{equation}
Next, using a similar notation, we see that
\begin{eqnarray*}
p_{x}^{(\theta)}(1) &=&\int_{0}^{1}P^{\theta}\bigl(O_{+}(x)\in dy\bigr)P^{\theta
}(T_{y}^{-}=\infty) \\
&=&c^{-}(\theta)\int_{0}^{1}P^{\theta}\bigl(O_{+}(x)\in dy\bigr)U_{-}^{\theta
}(y) \\
&=&c^{-}(\theta)\int_{0}^{1}P^{\theta}\bigl(O_{+}(x)>y\bigr)U_{-}^{\theta}(dy)
\\
&\leq&c^{-}(\theta)U_{+}^{\theta}(x)\int_{0}^{1}\overline{\mu}%
^{+}(\theta,y)U_{-}^{\theta}(dy) \\
:\!&=&c^{-}(\theta)I(\theta)U_{+}^{\theta}(x),
\end{eqnarray*}
where we have used the $P^{\theta}$ version of (\ref{bound}). Using
(\ref%
{er}) again gives the bound
%
%
\begin{equation} \label{one}
p_{x}^{(\theta)}(1)\leq xc_{0}c^{-}(\theta)I(\theta).
\end{equation}
Writing $O_{\pm}(n,x)$ for the successive overshoots upward and
downward over level~$x$, we then have
\[
p_{x}^{(\theta)}(n)=\int_{0}^{1}P^{\theta}\bigl(O_{-}(n-1,x)\in
dz\bigr)p_{z}^{(\theta)}(1)\leq c_{0}c^{-}(\theta)I(\theta)E^{\theta
}\bigl(O_{-}(n-1,x)\bigr).
\]
Also, Wald's identity gives $E^{\theta}(O^{-}(y))\leq m_{-}^{\theta
}U_{-}^{\theta}(y)$, where $m_{-}^{\theta}=E^{\theta}(H_{1}^{-})$,
and so
we have%
\begin{eqnarray*}
&&E^{\theta}\bigl(O_{-}(n-1,x)|O_{-}(n-2,x) =y\bigr)\\
&&\qquad=E^{\theta}(O^{-}(y))
\leq m_{-}(\theta)\int_{0}^{1}P^{\theta}(O_{y}^{+}\in
dz)U_{-}^{\theta
}(z) \\
&&\qquad=m_{-}(\theta)\int_{0}^{1}P^{\theta}(O_{y}^{+}>z)U_{-}^{\theta
}(dz) \\
&&\qquad\leq m_{-}(\theta)U_{+}^{\theta}(y)\int_{0}^{1}U_{-}^{\theta}(dz)%
\overline{\mu}^{+}(\theta,z) \\
&&\qquad\leq c_{0}m_{-}(\theta)I(\theta)y,
\end{eqnarray*}
where we have again used (\ref{bound}). Iterating this gives
%
%
\begin{equation} \label{e}
E^{\theta}\bigl(O_{-}(n-1,x)\bigr)\leq\{c_{1}(\theta)\}^{n-1}x,
\end{equation}
where $c_{1}(\theta)=c_{0}m_{-}(\theta)I(\theta)$, and thus%
\[
p_{x}^{(\theta)}(n)\leq c_{0}c^{-}(\theta)I(\theta)\{c_{1}(\theta
)\}^{n-1}x,\qquad n\geq1.
\]
Moreover, using (\ref{zero}) and (\ref{e}), we get the bound%
\begin{eqnarray*}
q_{x}^{(\theta)}(n) &=&\int_{0}^{1}P^{\theta}\bigl(O_{-}(n,x)\in
dz\bigr)P^{\theta
}(T_{z}^{+}=\infty) \\
&\leq&c_{0}k^{+}(\theta)E^{\theta}(O_{-}(n,x))\leq c_{0}k^{+}(\theta
)\{c_{1}(\theta)\}^{n-1}x.
\end{eqnarray*}
So, (\ref{lim}) will follow, provided that $\theta$ can be chosen such that
%
%
\begin{equation} \label{claim}
c_{1}(\theta)=c_{0}m_{-}(\theta)I(\theta)<1.
\end{equation}
To see this, we need to note first that $m_{-}(\theta)\leq
E(H_{1}^{-})$.
Also, provided that $k^{-}(\theta)\rightarrow\infty$, by applying
bound (\ref{er}) to $H^{-}$, we get $U_{-}^{\theta}(z)\rightarrow0$ for
each $%
z\in(0,1]$ as $\theta\rightarrow\infty$, and since $U_{-}^{\theta
}(z)\leq U_{-}(z)$ and $I<\infty$, dominated convergence will give
\[
I(\theta)=\int_{0}^{1}U_{-}^{\theta}(z)\mu^{+}(\theta,dz)\leq
\int_{0}^{1}U_{-}^{\theta}(z)\mu^{+}(dz)\rightarrow0\qquad\mbox{as }\theta
\rightarrow\infty.
\]
To see that $k^{-}(\theta)\rightarrow\infty$, note that the killing time
of $H^{-}$ under $P^{\theta}$ is the same as that of the ladder time
subordinator $L_{-}^{-1}$ and this has the distribution of
$L_{-}(\tau)$,
which is exp($\kappa_{-}(\theta))$, where $\kappa_{-}$ is the Laplace
exponent of $L_{-}$ under $P$. The assumption that $U_{-}(dx)>0$ for all
small $x>0$ implies that $L_{-}$ is not a compound Poisson process so, by
Corollary 3, page 17 of \cite{jb}, $\kappa_{-}(\infty)=\infty$ and, thus,
if we choose $\theta$ large enough, (\ref{claim}) will hold and the proof
is complete.
\end{pf}
\begin{proposition}
\label{B} \textup{(i)} Let $X$ be a L\'{e}vy process having $\Pi(\mathbb{R}%
^{+})=\infty$ and $\sigma>0$. Then a PRI exists.

\textup{(ii)} Let $X$ be a L\'{e}vy process having $\sigma=0$, $\Pi(\mathbb{R}%
^{+})=\infty$ and $\Pi(\mathbb{R}^{-})<\infty$. Then no PRI exists.
\end{proposition}
\begin{pf}
(i) Here, $\delta_{+}>0$ and $\delta_{-}>0$, so $U_{-}(x)\backsim
x/\delta
_{-}$ and since $\int_{0}^{1}x\mu^{+}(dx)$ is automatically finite, we
have $%
I<\infty$.

(ii) By the argument preceding Lemma \ref{O}, we can take $\Pi(\mathbb
{R}%
^{-})=0$ and assume that $\delta_{-}>0$, so that, again, $I$ is necessarily
finite. However, $\sigma=0$ and $\delta_{-}>0$ imply $\delta_{+}=0$,
so the
result follows.
\end{pf}

To deal with the remaining situations, we need the following lemma.
\begin{lemma}
\label{C} Let $X$ be an oscillating L\'{e}vy process whose L\'{e}vy measure
is supported by $[-1,1]$ and satisfies $\Pi([-1,0))=\Pi
((0,1]))=\infty$.
Suppose, additionally, that $\sigma=0$ and $\delta^{+}>0$. Then
$I=\int_{0}^{1}%
\overline{\mu}^{+}(x)U_{-}(dx)<\infty$ iff%
%
%
\begin{equation} \label{J2}
J=\int_{0}^{1}\frac{x^{2}\Pi(dx)}{ ( \int_{0}^{x}\int
_{y}^{1}\overline{%
\Pi}{}^{(-)}(s)\,ds\,dy ) ^{2}}<\infty.
\end{equation}
\end{lemma}
\begin{pf}
We use Vigons' ``\'{e}quation amicale invers\'{e}e'' (see \cite{vv}), which,
since our L\'{e}vy measure lives on $[-1;1]$, takes the form
\[
\overline{\mu}^{+}(x)=\int_{0}^{\infty}\overline{\Pi
}{}^{+}(x+y)U_{-}(dy)=%
\int_{x}^{1}U_{-}(y-x)\Pi(dy).
\]
We then use this in the following computation:%
\begin{eqnarray*}
I &=&\int_{0}^{1}\overline{\mu}^{+}(x)U_{-}(dx)<\infty
=\int_{0}^{1}\int_{x}^{1}U_{-}(y-x)\Pi(dy)U_{-}(dx) \\
&=&\int_{0}^{1}\Pi
(dy)\int_{0}^{y}U_{-}(y-x)U_{-}(dx)=\int_{0}^{1}U_{-}^{\ast2}(y)\Pi(dy).
\end{eqnarray*}
Next, we recall that the potential function $U_{-}(x)$ is increasing in $x$.
This is enough to show that%
\[
\bigl(U_{-}(y/2)\bigr)^{2}\leq U_{-}^{\ast2}(y)=\int
_{0}^{y}U_{-}(y-x)U_{-}(dx)\leq
(U_{-}(y))^{2}.
\]
Moreover, since $X$ oscillates, $H_{-}$ is an unkilled subordinator
with zero
drift and we have that $U_{-}(y)\approx y/A(y)$, where $A(y)=\int
_{0}^{y}%
\overline{\mu}^{(-)}(s)\,ds$ satisfies $A(y)/2\leq A(y/2)\leq A(y)$. This
implies that $U_{-}(y)\approx U_{-}(y/2)$ and therefore that
$U_{-}^{\ast
2}(y) \approx(U_{-}(y))^{2}$. We therefore conclude that%
\[
I=\int_{0}^{1}U_{-}^{\ast2}(y)\Pi(dy)<\infty\quad\Longleftrightarrow\quad
\int_{0}^{1}\frac{y^{2}\Pi(dy)}{A(y)^{2}}<\infty.
\]
Next, we need the ``\'{e}quation amicale int\'{e}gr\'{e}e'' of Vigon
(see \cite{vv}),
which, in our case, takes the form%
\[
\overline{\overline{\Pi}}{}^{(-)}(x)=\int_{x}^{1}\overline{\Pi}{}%
^{(-)}(y)\,dy=\int_{0}^{1}\overline{\mu}^{(+)}(y)\overline{\mu}%
^{(-)}(x+y)\,dy+\delta_{+}\overline{\mu}^{(-)}(x).
\]
Our assumptions imply that $\overline{\overline{\Pi}}{}^{(-)}(0+)>0$. If
$%
\overline{\overline{\Pi}}{}^{(-)}(0+)<\infty$, then it is obvious that
$0<%
\overline{\mu}^{(-)}(0+)<\infty$, and if $\overline{\overline{\Pi}}{}%
^{(-)}(0+)=\infty$, it is easy to deduce that $\overline{\mu}%
^{(-)}(0+)=\infty$. Then, from dominated convergence, it follows that%
\[
\lim_{x\downarrow0}\frac{\overline{\overline{\Pi}}{}^{(-)}(x)}{\overline
{\mu
}^{(-)}(x)}=\delta_{+}.
\]
Thus, in both cases,
$A(y)\approx\int_{0}^{y}\overline{\overline{\Pi}}{}
^{(-)}(z)\,dz$ and the result follows.
\end{pf}
\begin{pf*}{Proof of Theorem \ref{T}}
We have already covered all cases
except those having $\sigma=0$ and $\Pi(\mathbb{R}^{+})=\Pi(\mathbb
{R}%
^{-})=\infty$. By the\vspace*{1pt} standard argument, we can find another process,
$\tilde{%
X}$, which oscillates and whose L\'{e}vy measure $\tilde{\Pi}$ agrees
with $%
\Pi$ on $(-1,1)$ and is supported by\vspace*{1pt} $[-1,1]$, and is such that a PRI
exists for $X$ iff a PRI exists for $\tilde{X}$. Note that $\tilde{\Pi}%
([-1,0))=\tilde{\Pi}((0,1]))=\infty$ and that, in the obvious
notation, $%
\tilde{J}<\infty$ iff $J<\infty$. Proposition \ref{A} and Lemma \ref%
{C} then apply and show that a PRI exists iff $\delta_{+}>0$ and
$J<\infty$. If
$X$ has bounded variation, then $\overline{\overline{\Pi}}{}^{(-)}(0+)\in
(0,\infty)$, and $J=\infty$ is then automatic. If $X$ has unbounded
variation, then, as previously noted, $J<\infty$ implies $\delta
_{+}>0$ and
this completes the proof.
\end{pf*}
\begin{pf*}{Proof of Corollary \ref{Co}}
Since $\overline{\Pi
}{}^{-}(x)\approx
x^{-\alpha}$, where $1<\alpha<2$, we are in the unbounded variation case
and we need only check the value of the integral (\ref{J}). Clearly, $%
\int_{0}^{x}\int_{y}^{1}\overline{\Pi}{}^{(-)}(s)\,ds\,dy\approx
x^{2-\alpha}$,
so this reduces to checking whether%
\[
\int_{0}^{1}x^{2\alpha-2}\Pi(dx)=(2\alpha-2)\int_{0}^{1}x^{2\alpha
-3}%
\overline{\Pi}{}^{+}(x)\,dx<\infty
\]
and this holds iff $\beta<2\alpha-2$.
\end{pf*}
\begin{remark}
A similar calculation for the integral $L$ in (\ref{K}) shows that in this
example, $X$ creeps upward iff $\beta<\alpha$.
\end{remark}

\section{The excursion measure}

Evans \cite{se} and Winkel \cite{mw} both observed that we can
associate an
excursion theory with $K$.

They introduced $\Lambda_{t}=\inf\{x\dvtx K_{x}>t\}$, $Z=X-\Lambda$ and showed
that $Z$ is a strong Markov process with $\Lambda$ as a local time at zero.
It is clear that excursions away from 0 of $Z$ evolve in the same way as
excursions away from 0 of $X$, namely, they have the same semigroup,
but their
entrance laws will be different. For example, if $X=B$, then all
excursions of $Z$
are negative and the characteristic measure $n^{Z}$ is $n^{X}$
restricted to negative excursion paths.

Winkel showed that when $\sigma>0$, $n^{Z}$ is the restriction of $n^{X}$
to the set of excursion paths which \textit{start} negative. (To
do this, he had to demonstrate that all excursion paths either start negative
or start positive, that is, cannot leave $0$ in an oscillatory
fashion.) Therefore, $%
n^{Z}$ is absolutely continuous w.r.t. $n^{X}$.

However, this depends on both $\delta_{+}$ and $\delta_{-}$ being
positive. When $\sigma=0$ and $\delta_{+}>0$, we have $\delta_{-}=0$,
which means that excursions of $X$ have to return to 0 from below. By
time reversal, this means that they must start positive and since
excursions of $%
Z$ start negative, the two measures must be mutually singular whenever $
\sigma=0$. We believe that the problem of describing the excursion
measure $%
n^{Z}$ in this case is both interesting and difficult.

%

%
\printaddresses

\end{document}